\documentclass{elsart3-1}

\usepackage{amssymb}

\usepackage[english,francais]{babel}

\newtheorem{theorem}{Theorem}[section]

\newtheorem{e-proposition}[theorem]{Proposition}

\newtheorem{e-definition}[theorem]{Definition\rm}
\newtheorem{remark}{\it Remark\/}
\newtheorem{example}{\it Example\/}
\newtheorem{theoreme}{Th\'eor\`eme}[section]

\newtheorem{remarque}{\it Remarque}

\renewcommand{\theequation}{\arabic{equation}}
\def\og{\leavevmode\raise.3ex\hbox{$\scriptscriptstyle\langle\!\langle$~}}
\def\fg{\leavevmode\raise.3ex\hbox{~$\!\scriptscriptstyle\,\rangle\!\rangle$}}
\def\R{\mathbf{R}}
\def\C{\mathbf{C}}
\def\O{\mathcal{O}}
\def\e{\epsilon}

\newenvironment{prf}{\begin{trivlist}\item[]{\bf Proof:\ }}
{\mbox{}\hfill\rule{.08in}{.08in}\end{trivlist}}
\newenvironment{preuve}{\begin{trivlist}\item[]{\bf Preuve :\ }}
{\mbox{}\hfill\rule{.08in}{.08in}\end{trivlist}}
\journal{the Acad\'emie des sciences}
\begin{document}
% place in the next line the header (rubrique) chosen for your article,
% if you know it (you can also have 2, format : Header1/Header2
\centerline{Differential Geometry}
\begin{frontmatter}


%\title{Differentiability properties of the isoperimetric profile and topology of two dimensional analytic manifolds}
%\author{Renata Grimaldi\thanksref{label2}}
\selectlanguage{english}
\title{Differentiability of the isoperimetric profile and topology of analytic Riemannian manifolds}

% The [label1] can be suppressed if there is only one address for all authors

\selectlanguage{english}
\author[authorlabel1]{Renata Grimaldi},
\ead{grimaldi@unipa.it}
\author[authorlabel1]{Stefano Nardulli},
\ead{nardulli@unipa.it}
\author[authorlabel2]{Pierre Pansu}
\ead{Pierre.Pansu@math.u-psud.fr}
\ead[url]{http://www.math.u-psud.fr/$\sim$pansu}
\address[authorlabel1]{Dipartimento di Metodi e Modelli Matematici\\ Viale delle Scienze Edificio 8 - 90128 Palermo}
\address[authorlabel2]{Univ Paris-Sud, Laboratoire de
Math\'ematiques d'Orsay, Orsay, F-91405}
% If you know the dates of reception, and acceptation you can put them now;
%  idem the name of the person presenting the Note

\medskip
\begin{center}
{\small Received *****; accepted after revision +++++\\
Presented by ----}
\end{center}

\begin{abstract}
\selectlanguage{english}
% Text of abstract in English
We show that smooth isoperimetric profiles are exceptional for real analytic Riemannian manifolds. For instance, under some extra assumption, this can happen only on topological spheres. 
{\it To cite this article: Renata Grimaldi, Stefano Nardulli, Pierre Pansu, C. R. Acad. Sci. Paris, Ser. I *** (200+).}

\vskip 0.5\baselineskip

\selectlanguage{francais}
% Text of abstract in French

\noindent{\bf R\'esum\'e} \vskip 0.5\baselineskip \noindent

{\bf Differentiabilit\'e du profil isop\'erim\'etrique et topologie des vari\'et\'es riemanniennes analytiques r\'eelles}
On montre que la differentiabilit\'e du profil isop\'erim\'etrique est une condition tr\`es contraignante pour les vari\'et\'es riemannniennes analytiques r\'eelles. Par exemple, sous une hypoth\`ese suppl\'ementaire, ce n'est possible que si la vari\'et\'e est hom\'eomorphe \`a une sph\`ere. 
{\it Pour citer cet article : Renata Grimaldi, Stefano Nardulli, Pierre Pansu, C. R. Acad. Sci. Paris, Ser. I *** (200+).}

\end{abstract}
\end{frontmatter}

% now the Version franaise abr\'eg\'ee, if it exists
\selectlanguage{francais}
\section*{Version fran\c{c}aise abr\'eg\'ee}
% Text of your Version franaise abr\'eg\'ee here.
% Note you do not need to repeat here equations that you use in the
% main text - for example 'voir (3)' is quite acceptable.

\selectlanguage{francais}

Le {\em profil isop\'erim\'etrique} d'une vari\'et\'e riemannienne $M$ est la fonction $I$ qui, en $v\in (0,\mathrm{vol}(M))$, vaut la borne inf\'erieure $I(v)$ des volumes des bords des domaines de volume \'egal \`a $v$.

Il y a peu de vari\'et\'es riemanniennes dont on connaisse exactement le profil (voir \cite{B} pour un \'etat de l'art). Parmi tous les exemples connus, seule la sph\`ere ronde a un profil lisse. Toute vari\'et\'e poss\`ede-t-elle une m\'etrique riemanienne de profil lisse ? Dans \cite{GP}, on construit sur toute vari\'et\'e de dimension 2 ayant suffisamment de bout des m\'etriques lisses de profil lisse. Il est probable que des exemples similaires existent sur de nombreuses vari\'et\'e compactes. En revanche, il semble plus difficile de construire des exemples analytiques r\'eels.

On sait (\cite{D}) que le maximum de la courbure scalaire intervient dans le d\'eveloppement asymptotique du profil pr\`es de 0. Une condition suffisante pour que le profil soit lisse au voisinage de 0 est que ce maximum soit atteint en un seul point, et que celui-ci soit un point critique non d\'eg\'en\'er\'e, voir \cite{N}. Sous cette hypoth\`ese suppl\'ementaire, nous montrons que demander un profil lisse est tr\`es restrictif.

\begin{theoreme}\label{main}Soit $M$ une vari\'et\'e riemannienne compacte analytique r\'eelle. On fait les hypoth\`eses suivantes.
\begin{enumerate}
  \item Toutes les solutions du probl\`eme isoperim\'etrique dans $M$ sont lisses.      
  \item La courbure scalaire atteint son maximum en un seul point, et celui-ci est un point critique non d\'eg\'en\'er\'e.
  \item Le profil isop\'erim\'etrique de $M$ est lisse.
\end{enumerate}
Alors $M$ est hom\'eomorphe \`a une sph\`ere.
\end{theoreme}

\begin{remarque}
L'hypoth\`ese (i) est toujours satisfaite si $\mathrm{dim}(M)<8$. Nous ignorons si on peut se passer de l'hypoth\`ese (ii).
\end{remarque}

\begin{preuve}
Si le profil est lisse, il est analytique. Les solutions lisses du probl\`eme isop\'erim\'etrique forment alors un espace analytique r\'eel de dimension finie. Par hypoth\`ese, les solutions de petit volume forment un unique arc, une famille \`a un param\`etre (on peut prendre le volume comme param\`etre) de domaines proches des petites boules centr\'ees au maximum de la courbure scalaire. Cet arc analytique se prolonge analytiquement, et ne peut se terminer que par la famille des domaines dont les pr\'ec\'edents sont les compl\'ementaires. $M$ est donc partag\'ee en deux domaines diff\'eomorphes \`a des boules ferm\'ees.
\end{preuve}

\selectlanguage{english}
% main text
\section{Introduction}
\label{intro}

The {\em isoperimetric profile} of a Riemannian manifold $M$ is the function $I$ which maps $v\in (0,\mathrm{vol}(M))$ to the infimum $I(v)$ of boundary volumes of domains of volume equal to $v$.

Very few profiles of compact Riemannian manifolds are exactly known (\cite{B} contains a state of the art). Among those, only the round sphere has a smooth profile. We address here the differentiability question: does every compact manifold admit a metric with smooth profile ? In \cite{GP}, smooth metrics with smooth profiles are constructed on 2-dimensional manifolds, provided they have enough ends. It is likely that similar examples can be found on many compact manifolds. However, it seems much harder to produce real analytic metrics with smooth profiles.

It is known (\cite{D}) that the asymptotics of the isoperimetric profile near 0 is determined by the maximum of scalar curvature. A sufficient condition for the profile to be smooth near 0 is that this maximum be achieved at only one point, which is a nondegenerate critical point, see \cite{N}. Under this assumption, we show that having a smooth profile is very restrictive.

\begin{theorem}\label{main}Let $M$ be a compact real analytic Riemannian manifold. Assume the following.
\begin{enumerate}
  \item All solutions of the isoperimetric problem in $M$ are smooth.
  \item The scalar curvature function achieves its maximum at a unique point, and this point is a nondegenerate critical point.
  \item The isoperimetric profile of $M$ is smooth.
\end{enumerate}
Then $M$ is homeomorphic to a sphere.
\end{theorem}

\begin{remark}
Assumption (i) is always satisfied if $\mathrm{dim}(M)<8$. We do not know wether assumption (ii) can be removed.
\end{remark}

% etc, etc

% \appendix
\section{Real analytic spaces}

We collect more or less standard material on real analytic spaces. We express thanks to Jean-Jacques Risler for his help with these matters.

\subsection{Definitions}

The following definitions are borrowed from H. Hironaka, \cite{H}, sections 1 and 5.

\begin{e-definition}
\label{ringed}
{\em (H. Hironaka, \cite{H}, Definition 1.1).} A {\em ringed space} is a topological space equipped with a sheaf of commutative algebras over $\R$. A {\em morphism} of ringed spaces $f:(X,\O_X)\to(Y,\O_Y)$ is the data of a continuous map $f:X\to Y$ together with a homomorphism of sheaves $f_* \O_X \to\O_Y$. The stalk $\O_{X,x}$ at each point is a local ring. Its Krull dimension (the length of the longest chain of prime ideals in the maximal ideal) is called the {\em dimension of $X$ at $x$}.
\end{e-definition}

\begin{example}
\label{exringed}
Let $U\subset\R^n$ be an open set, let $g:U\to\R^m$ be a real analytic map. Then $U\cap g^{-1}(0)$ equipped with the sheaf of restrictions (to open subsets of $U\cap g^{-1}(0)$) of real analytic functions on open subsets of $\R^n$, is a ringed space. Such ringed spaces are called {\em local models}.
\end{example}

\begin{e-definition}
\label{defsmooth}
{\em (H. Hironaka, \cite{H}, Definition 5.6).} A local model is {\em smooth} at a point $x$ if the defining map $g$ is a submersion at $x$.
\end{e-definition}

\begin{e-definition}
\label{raspace}
{\em (H. Hironaka, \cite{H}, Definition 1.5).} A {\em real analytic space} is a ringed space in which every point has a neighborhood isomorphic to some local model. Its dimension is the supremum of dimensions at points. Smooth points are those where the local model is smooth.
\end{e-definition}

%\begin{prop}
%\label{filtration}
%{\em (H. Hironaka, \cite{H}, Proposition 5.8).} Every $\sigma$-compact real analytic space admits a locally finite filtration $X=X_0 \supset X_1 \supset\cdots$ by closed real analytic subspaces $X_i$ such that $X_i \setminus X_{i+1}$ is smooth. Furthermore, for $x\in X_{i+1}$, $\mathrm{dim}_x (X_{i+1})<\mathrm{dim}_x (X_{i})$.
%\end{prop}

\subsection{Extracting a 1-dimensional subset}

\begin{prop}
\label{filtrationbis}
Let $X$ be an $n$-dimensional paracompact real analytic space. Then $X$  admits a locally finite filtration $X=X_n \supset X_{n-1} \supset\cdots\supset X_0$ by closed real analytic subspaces $X_i$ such that $X_i \setminus X_{i-1}$ is smooth of dimension $i$.
\end{prop}

\begin{prf}
The proof is an adaptation of Proposition 5.8 of \cite{H}. $X$ admits a complexification $X^{\C}$. This is a complex analytic space containing $X$ equipped with an antiholomorphic involution $\sigma$ whose fixed point set equals $X$. The set $X_{n-1}^{\C}$ of points of $X^{\C}$ where $X^{\C}$ is not smooth and $n$-dimensional is defined in local models by holomorphic equations. Therefore $X_{n-1}^{\C}$ is a complex analytic subspace. Let $X_{n-1}=X\cap X_{n-1}^{\C}$. Then $X\setminus X_{n-1}$ is smooth of dimension $n$. Define other $X_i$'s recursively.
\end{prf}

\begin{cor}
\label{extract}
Let $X$ be a paracompact real analytic space. Let $x$ be a smooth point where $X$ has dimension 1. There exists a closed real analytic subset $Y\subset X$ which contains a neighborhood of $x$ and whose dimension at each point is 1.
\end{cor}

\begin{prf}
Put $Y=X_1$. By construction of $X_{n-1},\ldots,X_1$, the open set $U$ of points where $X$ is smooth of dimension 1 is contained in $Y$.
\end{prf}

\subsection{Endpoints}

\begin{prop}
\label{even}
{\em (D. Sullivan, \cite{S}, page 166, example a)).} Let $X$ be a 1-dimensional real analytic space. Let $x\in X$ be a singular point. Then a punctured neighborhood of $x$ in $X$ consists of an even number of arcs. They match pairwise into real analytic arcs.
\end{prop}

\begin{cor}
\label{ends}
Let $X$ be a real analytic space whose dimension at each point is 1. Assume that the singular set of $X$ is compact. Then each connected component of $X$ has a finite, even number of ends.
\end{cor}

\begin{prf}
Let $S\subset X$ denote the singular set of $X$. Since $S$ is a closed real analytic subspace, $S$ is discrete, and therefore finite. Since the property of being a real analytic space is local, and local models are locally connected, one can assume that $X$ is connected. $X\setminus S$ is a 1-dimensional manifold, thus a collection of arcs and circles. Since $X$ is connected, either $S$ is empty, $X$ is a circle (no ends) or there are no circles. Then each arc has at least one endpoints in $S$, thus there are finitely many arcs. Ends correspond to ends of arcs which do not belong to $S$. Let $N$ denote the number of ends of $X$ and $P$ the number of pairs $(s,a)$ where $s\in S$ and $a$ is an arc with $s$ as an endpoint. Then $N+P$ is twice the number of arcs, $P$ is even (according to Proposition \ref{even}), thus $N$ is even. 
\end{prf}

\section{Constant mean curvature hypersurfaces}

\begin{e-definition}
\label{top}
Fix $\alpha\in(0,1)$. 
Let us define the $C^{2,\alpha}$-topology on the space of domains with smooth boundary as follows: as neighborhoods of a smooth domain $B$, take all domains $S$ whose boundary is the graph, in normal exponential coordinates, of a $C^{2,\alpha}$-small function on $\partial B$. 
\end{e-definition}

\begin{e-proposition}
\label{mean}
Let $M$ be a compact real analytic Riemannian manifold. Let $f$ be a positive real analytic function on the interval $(0,\mathrm{vol}(M))$. The set $X^f$ of pairs $(B,h)$ where $B$ is a smooth domain whose boundary has constant mean curvature equal to $h$ and satisfying $\mathrm{vol}(\partial B)=f(\mathrm{vol}(B))$ has a natural structure of a finite dimensional real analytic space.
\end{e-proposition}

\begin{prf}
The proof is essentially contained in section 4 of \cite{GNP}. 
Let $(B,k)$ belong to $X^f$. By definition, a neighborhood of $B$ consists of domains $D_u$ whose boundaries are graphs in normal coordinates of $C^{2,\alpha}$ functions $u$ which satisfy the constant mean curvature equation $H(u)=h$, $h\in\R$. Since the linearized operator at $u=0$, $L_{B}$, need not be invertible, one solves instead
\begin{eqnarray*}
\Phi_{B}(u,h)=P_B (H_{B}(u)-h)=0,
\end{eqnarray*}
where $P_B$ is the orthogonal projection onto the $L^2$-orthogonal complement $K_{B}^{\bot}$ of the kernel of $L_{B}$ in $C^{0,\alpha}(\partial B)$. Then 
\begin{eqnarray*}
\Phi_{B}:C^{2,\alpha}(\partial B)\times\R\to K_{B}^{\bot}
\end{eqnarray*}
is a real analytic map, whose linearization at 0 is $P_B \circ L_{B}$. By construction, it is onto with a finite dimensional kernel. Therefore the set of solutions is a finite dimensional smooth real analytic submanifold $Y_{B}$ of $C^{2,\alpha}(\partial B)\times\R$ in a neighborhood of $B$. On $Y_{B}$, 
\begin{eqnarray*}
(u,h)\mapsto (H(u)-h,\mathrm{vol}(\partial D_u)-f(\mathrm{vol}(D_u)))
\end{eqnarray*}
is a real analytic map, whose zero set, denoted by $X_{B}$, is a neighborhood of $B$ in $X^f$. A choice of coordinates for $Y_{B}$ maps $X_{B}$ to a local model. The sheaf $\O_{X^f}$ is defined on small enough open sets by restricting to $X_{B}$ analytic functions on $C^{2,\alpha}(\partial B)\times\R$. Different choices of $B$ yield local embeddings of $X^f$ into Banach spaces which differ only up to a change of chart for the (infinite dimensional) manifold structure on the space of domains. Since constant mean curvature hypersurfaces in real analytic Riemannian manifolds are real analytic, these changes of charts are real analytic diffeomorphisms between open sets of $C^{2,\alpha}$ spaces. This yields the structure of a real analytic space on $X^f$.
\end{prf}

\section{Compactness}

\begin{e-proposition}\label{compact}Let $M$ be a compact Riemannian manifold with isoperimetric profile $I$. Assume that all solutions of the isoperimetric problem in $M$ are smooth. Then for each $\e>0$, the subset
\begin{eqnarray*}
X^{I,\e}=\{B\in X^I \,|\,\e\leq \mathrm{vol}(B)\leq\mathrm{vol}(M)-\e\}
\end{eqnarray*}
is compact in $C^{2,\alpha}$ topology.
\end{e-proposition}

\begin{prf} The proof is essentially the same as in lemma $8$ of \cite{GNP}. It depends on the fact that the $C^{2,\alpha}$ topology is equivalent to the flat topology of Geometric Measure Theory on the space of smooth isoperimetric domains. Compactness follows the compactness theorem for currents of bounded mass and boundary mass, provided one assumes that flat limits of smooth minimizers (they are again minimizers) are smooth.
\end{prf}

\section{Proof of Theorem \ref{main}}

\begin{prf} 
By definition, elements of $X^I$ coincide with solutions of the isoperimetric problem. We know from \cite{N} that elements of $X^I$ of small volume are nearly round balls emanating from the unique maximum of scalar curvature (see also \cite{Y}). They form a smooth connected 1-dimensional manifold. In other words, for $\e$ small enough, $X^I$ is the union of $X^{I,\e}$ and of two smooth one-ended arcs, 
\begin{eqnarray*}
X^{I,\leq\e}=\{B\in X^I \,|\,\mathrm{vol}(B)\leq\e\}\quad \textrm{and}\quad X^{I,\geq V-\e}=\{B\in X^I \,|\,\mathrm{vol}(B)\geq \mathrm{vol}(M)-\e\}.
\end{eqnarray*}
>From Proposition \ref{compact}, it follows that $X^I$ has exactly two ends.

According to \cite{GNP}, $I$ is semi-analytic on $[0,\mathrm{vol}(M)]$. If smooth on $(0,\mathrm{vol}(M))$, it must be analytic on the same interval. Consider the subset $D\subset X^I$ consisting of domains diffeomorphic to the $n$-ball. Since $D$ is open and closed in $X^I$, it inherits the structure of a finite dimensional real analytic space. Let $D=D_N \supset D_{N-1}\supset\cdots\supset D_1 \supset D_0$ be the filtration provided by Proposition \ref{filtrationbis}. The subset $D_1$ is a real analytic space of dimension 1. Let $D'_1$ denote the connected component of $D_1$ which contains $X^{I,\leq\e}$.

By corollary \ref{ends}, $D'_1$ has a finite even number of ends. Since it contains one of the two ends of $X^I$, it must contain the other. In other words, there exists a domain $B$ in $M$ which is diffeomorphic to an $n$-ball, and whose complement is diffeomorphic to an $n$-ball as well. This implies that $M$ is homeomorphic to a sphere. 
\end{prf}

% The Acknowledgements are an un-numbered section
%\section*{Acknowledgements}
% Acknowledgements text here

\end{document}